\documentclass[12pt]{article}
\usepackage[pctex32]{graphics}
\usepackage{amsthm,amsmath,amssymb,amscd,verbatim,epsfig}
\usepackage{graphics,pifont}
\usepackage{amssymb}         
\usepackage{graphicx}
\newcommand{\beq}{\begin{equation}}
\newcommand{\eeq}{\end{equation}}

\date{}
\newcommand{\al}{\alpha}
\newcommand{\ga}{\gamma}
\newcommand{\be}{\beta}

\newcommand{\si}{\sigma}
\newcommand{\Si}{\Sigma}

\newcommand{\dl}{\delta}
\newcommand{\emp}{\emptyset}
\newcommand{\iy}{\infty}

\newcommand{\lra}{\longrightarrow}
\newcommand{\sub}{\subset}
\newcommand{\supp}{\supset}

\newcommand{\sq}{$\blacksquare$}

\begin{document}

\title{On\ the\ Invariance\ of\ Li-Yorke\ Chaos\ of\ Interval\ Maps}
\author{Bau-Sen Du \\ [.3cm]
Institute of Mathematics \\
Academia Sinica \\
Taipei 11529, Taiwan \\
dubs@math.sinica.edu.tw \\}
\maketitle

\begin{abstract}
In their celebrated "Period three implies chaos" paper, Li and Yorke proved 
that if a continuous interval map $f$ has a period 3 point then there is an uncountable
scrambled set $S$ on which $f$ has very complicated dynamics.  One question
arises naturally: Can this set $S$ be chosen invariant under $f$?  The answer is
positive for turbulent maps and negative otherwise.  In this note, we shall
use symbolic dynamics to achieve our goal.  In particular, we 
obtain that the tent map $T(x) = 1 - |2x-1|$ on $[0, 1]$ has a dense uncountable 
invariant 1-scrambled set which consists of transitive points.
\end{abstract}


\section{Introduction.}  Let $I$ be a real compact interval and let $f : I \lra I$ be a continuous map.  For each integer $n \ge 1$, let $f^n$ be defined by: $f^1 = f$ and $f^n = f \circ f^{n-1}$ when $n \ge 2$.  For every $c$ in $I$, we call the set $O_f(c)=\{ \, f^k(c) \mid k \ge 0 \, \}$ the {\it{orbit}} of $c$ (under $f)$ and call $c$ a {\it{periodic point}} of $f$ with least period $m$ if $f^m(c) = c$ and $f^i(c) \ne c$ whenever $0 < i < m$.  It is well known {\bf{\cite{bc}}}, {\bf{\cite{du1}}}, {\bf{\cite{ly}}} that if $f$ has a periodic point of least period not a power of 2, then there exist a number $\dl > 0$ and an uncountable (called \, $\dl$-scrambled) set $S$ of $f$ such that (i) for any $x \ne y$ in $S$, $\limsup_{n \to \iy} |f^n(x)-f^n(y)| \ge \dl$ and $\liminf_{n \to \iy} |f^n(x)-f^n(y)| =0$; and (ii) for any $x$ in $S$ and any periodic point $p$ of $f$, $\limsup_{n \to \iy} |f^n(x)-f^n(p)| \ge \dl/2$.  One question arises naturally: Can $S$ be chosen invariant under $f$?  The answer is negative.  Indeed, for any continuous map $f$ from $[0, 1]$ into itself, we define another continuous map $F$ from $[0, 3]$ into itself by letting $F(x) = f(x) + 2$ for $0 \le x \le 1$; $F(x) = (2-x)[f(1)+2]$ for $1 \le x \le 2$; and $F(x) = x-2$ for $2 \le x \le 3$.  This function $F$ is generally called the {\it{square root of $f$}} and the process from $f$ to $F$ is called the {\it{square root construction}}.  If $F$ has a non-empty invariant scrambled set $S$, then since every non-fixed point of $(1, 2)$ will spiral out and eventually enter $[0, 1] \cup [2, 3]$ and since $F([0, 1]) \sub [2, 3]$ and $F([2, 3]) = [0, 1]$, we see that $S \cap [0, 1] \ne \emp$ and $S \cap [2,3] \ne \emp$.  Choose any $x$ in $S \cap [0, 1]$ and any $y$ in $S \cap [2, 3]$.  Since, for each $n \ge 1$, $f^n(x)$ and $f^n(y)$ lie on opposite sides of $[1, 2]$, we obtain that $\liminf_{n \to \iy} |f^n(x) - f^n(y)| \ge 1$.  Therefore, $F$ cannot have a non-empty invariant scrambled set.  In this note, we shall use symbolic dynamics to show that if $f$ is turbulent then $f$ has an uncountable invariant scrambled set in the recurrent set.  Some related problems are also discussed.

\section{Turbulent maps have uncountable invariant scrambled sets.}
We say that $f$ is {\it{turbulent}} [{\bf 1}, p.25] if there exist two closed subintervals $J_0$ and $J_1$ of $I$ with at most one common point such that $f(J_0) \cap f(J_1) \supp J_0 \cup J_1$.  We say that $f$ is {\it{strictly turbulent}} if $J_0$ and $J_1$ can be chosen disjoint.  Let $f$ be a turbulent map with $J_0$ and $J_1$ as stated above and $\max J_0 \le \min J_1$.  Since $f(J_0) \supp J_0 \cup J_1$, there exist two points $u$ and $v$ in $J_0$ such that $f(u) = \min J_0$ and $f(v) = \max J_1$.  Without loss of generality, suppose $u < v$.  Let $z = \max \{ \, x \in [u, v] \, \mid \, f(x) = x \, \}$, $b = \min \{ \, v \le x \le \max J_1 \, \mid \, f(x) = z \, \}$,  $c = \min \{ \, z < x \le v \, \mid \, f(x) = b \, \}$ and $a = \max \{ \, v \le x \le b \, \mid \, f(x)= b \, \}$.  Then $z < c \le a < b$ and we obtain [{\bf 1}, p.26] that (i) $f(z) = z = f(b)$ and $f(c) = b = f(a)$; (ii) $x < f(x) < b$ for $z < x < c$; and (iii) $z < f(x) < b$ for $a < x < b$.  Let $I(0) = [z, c]$ and $I(1) = [a, b]$.  Then $I(0)$ and $I(1)$ have at most one point in common and $f(I(0)) \cap f(I(1)) \supp I(0) \cup I(1)$.

For $\al_i = 0$ or 1, let $I(\al_0\al_1 \cdots \al_n)$ be any closed subinterval of $I(\al_0\al_1 \cdots \al_{n-1})$ of minimum length {\bf{\cite{bc}}}, {\bf{\cite{ev}}} such that $f(I(\al_0\al_1 \cdots \al_n)) = I(\al_1\al_2 \cdots \al_n)$.  Hence, $f$ maps the endpoints of $I(\al_0\al_1 \cdots \al_n)$ onto those of $I(\al_1\al_2 \cdots \al_n)$ and maps the interior of $I(\al_0\al_1\cdots \al_n)$ onto the interior of $I(\al_1\al_2 \cdots \al_n)$.  Let $\Sigma_2 = \{ \, \al \, \mid \, \al =  \al_0\al_1 \cdots$, where $\al _i = 0$ or 1 $\}$ be the metric space with metric $d$ defined by $d(\al_0\al_1 \cdots$, $\be_0\be_1 \cdots) = \sum_{i=0}^\iy {|\al_i-\be_i|}/{2^{i+1}}$ and let $\si$ be the shift map defined by $\si(\al_0\al_1 \cdots) = \al_1\al_2 \cdots$.  For any $\al = \al_0\al_1 \cdots$ in $\Si_2$, let $I(\al) = \bigcap_{n=0}^{\iy} I(\al_0\al_1 \cdots \al_n)$.  Then each $I(\al) \sub I(\al_0)$ is either a compact interval or consists of one point.  Furthermore, $f(I(\al)) = I(\si \al)$ and $f$ maps the endpoints of $I(\al)$ onto those of $I(\si \al)$ if both $I(\al)$ and $I(\si \al)$ are nondegenerate intervals.  So, if $\al \in \Si_2$ is a periodic point of $\si$ with least period $m$ and if $I(\al) = [x, y]$ ($x$ may equal $y$), then both $x$ and $y$ are periodic points of $f$ with least period $m$ or $2m$.  Now let $\bar 0 \in \Si_2$ denote the sequence consisting of all $0$'s.  Then $f(I(\bar 0)) = I(\bar 0) \sub I(0)$.  Since $z$ is the unique fixed point of $f$ in $I(0)$ and $z = \min I(0)$, $I(\bar 0) = [z, w]$ for some $w$ in $I(0)$.  Since $x < f(x)$ for all $z < x$ in $I(0)$, we actually have $I(\bar 0) = \{ \, z \, \}$.  Since $f(I(1 \bar 0)) = I(\bar 0) = \{ \, z \, \}$ and since $b$ is the only point in $I(1)$ mapping to $z$, we obtain that $I(1\bar 0) = \{ \, b \, \}$.   Note that, if $c < a$, then $I(\al)$ and $I(\be)$ are disjoint for any $\al \ne \be$ in $\Si_2$.  If $c = a$ and $I(\al) \cap I(\be) \ne \emp$ for some $\al \ne \be$ in $\Si_2$, then for some $k \ge 0$ and $\ga_i = 0$ or 1, $0 \le i \le k-1$, we have $\{ \al, \be \} = \{ \ga_0\ga_1 \cdots \ga_{k-1}01\bar 0, \ga_0\ga_1 \cdots \ga_{k-1}11\bar 0 \}$ {\bf [5}, p.274].

Let $\al$ be a point of $\Si_2$ whose expansion $\al_0\al_1\cdots $ contains every finite sequence of $0$'s and $1$'s.  Then $\al$ has a dense orbit in $\Si_2$ (points with dense orbits are called {\it{transitive points}}).  For any integer $k \ge 1$, let $0^k \,\, (1^k$ \, respectively) denote the finite sequence of $k$ consecutive $0$'s ($1$'s respectively).  For any $\be = \be_0\be_1 \cdots$ in $\Si_2$, let $$\omega_\be = \al_0\be_0 \, 1 \,\,\, \al_0\al_1\be_00\be_10 \, 10 \,\,\, \cdots \,\,\, \al_0\al_1 \cdots \al_k \be_00^k\be_10^k \cdots \be_k0^k \, 1 0^k \, \cdots.$$  Then $f^{m(m+1)(2m+7)/6}(I(\omega_\be)) \sub I(\al_0\al_1 \cdots \al_m \be_00^m\be_10^m \cdots \be_m0^m \, 1 0^m)$ for any $m \ge 1$.  Since there can have at most countably many $\eta = \eta_0\eta_1 \cdots$ in $\Si_2$ such that $I(\eta)$ is an interval, there exists an uncountable subset $B$ of  $\Si_2$ such that if $\be$ is in $B$ then \,$I(\omega_\be)$ consists of exactly one point, say, $I(\omega_\be) = \{ \, x_\be \, \}$.  Let $W$ be the collection of such points $x_\be$ for all $\be$ in $B$ and let $S = \bigcup_{i \ge 0} f^i(W)$.  Then $S$ is an uncountable invariant set of $f$.    

For any $m \ge 1$, let $s_m = m(m+1)(2m+7)/6$.  For any $\be, \ga$, and $\eta$ in $\Si_2$ with $\be \ne \ga$, if $\be_k \ne \ga_k$ for some $k \ge 0$, then for any $i \ge 0$ and all $m > k$, $f^{s_m + (k+1)(m+1) - i}(f^i(x_\be)) \in I(\be_k0^m)$ and $f^{s_m + (k+1)(m+1) - i}(f^i(x_\ga)) \in I(\ga_k0^m)$.  Therefore, $\limsup_{n \to \iy} |f^n(f^i(x_\be)) - f^n(f^i(x_\ga))|$  $\ge \text{dist} (I(\bar 0), I(1 \bar 0)) =|z-b|> 0$.  Furthermore, for any $0 \le i < j$ and all $m > j$, $f^{s_{m+1}-m-1-i}(f^i(x_\be)) \in I(10^m)$ and $f^{s_{m+1}-m-1-i}(f^j(x_\eta)) \in I(0^{m-j+i+1})$.  Thus, $\limsup_{n \to \iy} |f^n(f^i(x_\be)) - f^n(f^j(x_\eta))| \ge |z-b|$.  This shows that, for any $x \ne y$ in $S$, $\limsup_{n \to \iy} |f^n(x) - f^n(y)| \ge |z-b|$.  On the other hand, we have $f^{s_{m+1}-m-i}(f^i(x_\be)) \in I(0^m) \sub I(0^{m+i-j})$, $f^{s_{m+1}-m-i}(f^i(x_\ga)) \in I(0^m) \sub I(0^{m+i-j})$, and $f^{s_{m+1}-m-i}(f^j(x_\eta))$ $\in I(0^{m+i-j})$.  Since the sequence $<I(0^{m+i-j})>$ shrinks to the fixed point $z$, we obtain that  $\liminf_{n \to \iy}$ $|f^n(f^i(x_\be)) - f^n(f^i(x_\ga))| = 0$ and $\liminf_{n \to \iy} |f^n(f^i(x_\be)) - f^n(f^j(x_\eta))| = 0$.  This shows that, for any $x \ne y$ in $S$, $\liminf_{n \to \iy} |f^n(x) - f^n(y)| = 0$.  Since, for any $x$ in $S$, the iterates $f^n(x)$ approach $b$ infinitely often and stay as close to the fixed point $z$ of $f$ for as long as we wish, we obtain that, for any periodic point $p$ of $f$, $\limsup_{n \to \iy}$ $|f^n(x) - f^n(p)| \ge |z-b|/2$.  This proves that $S$ is an uncountable invariant $|b-z|$-scrambled set of $f$.

Since, for every $\be$ in $B$, $\omega_\be$ is a transitive point of $\Si_2$ and so, the iterates $\si^n(\omega_\be)$ approach every $\ga$ in $\Si_2$ infinitely many times.  Consequently, the iterates $f^n(x_\be)$ of the unique point $x_\be$ of $I(\omega_\be)$ approach the set $I(\ga)$ infinitely many times.  Therefore, the $\omega$-limit set $\omega(x_\be, f)$ contains $I(1\bar 0) \cup I(\bar 0) = \{b, z \}$ and contains the point $x_\ga$ for all $\ga \in B$ and at least one endpoint of $I(\omega_\ga)$ for all $\ga \notin B$.  In particular, $x_\be \in \omega(x_\be, f)$ (such points $x_\be$ are called {\it{recurrent points}}).  Therefore, we have obtained the following result.

\noindent
{\bf Theorem 1.}
{\it Let $f$ be a turbulent map.  Then there exist a fixed point $z$ of $f$, a point $b \ne z$ with $f(b) = z$, and an uncountable invariant $|b-z|$-scrambled set $S \sub [z : b]$ in the recurrent set of $f$, where $[z:b]$ denotes the closed interval with $z$ and $b$ as endpoints, such that the set $\bigcap_{x \in S} \omega(x, f)$ contains $\{ z, b \} \cup S$ and contains a periodic point of least period $n$ or $2n$, for every positive integer $n$.}

By mimicking the proof of Theorem 1, we obtain the following result {\bf{\cite{de}}}.

\noindent
{\bf Theorem 2.}
{\it Let $0 < a < 1$ and let $f$ be a continuous map from $[0, 1]$ onto itself such that (i) $f(0) = 0 = f(1)$ and $f(a) = 1$ and (ii) $f$ is strictly increasing on $[0, a]$ and stricly decreasing on $[a, 1]$.  Then the following statements are equivalent:
\begin{itemize}
\item[\rm{(a)}]
$f$ has a dense orbit.

\item[\rm{(b)}]
$f$ has dense periodic points.

\item[\rm{(c)}]
$f$ has sensitive dependence on initial conditions; i.e., there exists a positive number $\dl$ such that for any point $x$ in $[0, 1]$ and any open neighborhood $V$ of $x$ there exist a point $y$ in $V$ and a positive integer $n$ such that $|f^n(x)-f^n(y)| \ge \dl$.
\end{itemize}

\noindent
Furthermore, if $f$ has a dense orbit, then $f$ is topologically conjugate to the tent map $T(x) = 1 - |2x-1|$ on $[0, 1]$ and $f$ has a dense uncountable invariant 1-scrambled set which consists of transitive points.}

\section{Extensions to some non-turbulent maps.}  
Since the recurrent set of $f^m$ is identical to that of $f$ {\bf [1}, p.82], we can, by mimicking the proof of Theorem 1, extend the result of Theorem 1 to some non-turbulent maps.

\noindent
{\bf Theorem 3.}
{\it Let $g = f^m$, where $m \ge 2$.  Assume that $g$ is turbulent and let $z < c \le a < b$ be points such that (i) $g(z) = z = g(b)$ and $g(c) = b = g(a)$; (ii) $x < g(x) < b$ for $z < x < c$; and (iii) $z < g(x) < b$ for $a \le x < b$.  If $z$ is also a fixed point of $f$, then for some $\dl > 0$, $f$ has an uncountable invariant $\dl$-scrambled set $W$ in the recurrent set of $f$ such that the set $\bigcap_{x \in W} \omega(x, f)$ contains $W$ and contains periodic points of $f$ of arbitrarily large periods.}

\noindent
{\it Proof.}
Let $S \sub [z, b]$ be an uncountable invariant $|b-z|$-scrambled set of $g$ obtained as described in the proof of Theorem 1 which is in the recurrent set of $g$.  Let $W = \cup_{j=0}^{m-1} f^j(S)$.  Let $x$ and $y$ be any two points in $S$ and let $0 \le s < t \le m-1$ be integers.  Let $< k_i >$ be an increasing sequence of positive integers.  If $\lim_{i \to \iy} g^{k_i}(x) = b$ and $\lim_{i \to \iy} g^{k_i}(y) = z$, then since $z$ is also a fixed point of $f$, $\lim_{i \to \iy} f^{k_i+t-s}(y) = z$.  So, $\lim_{n \to \iy} |f^n(f^s(x) - f^n(f^t(y)|$ $\ge \lim_{i \to \iy} |f^{mk_i}(x) - f^{mk_i+t-s}(y)| = |b-z|$.  If $\lim_{i \to \iy} g^{k_i}(x) = z$ and $\lim_{i \to \iy} g^{k_i}(y) = b$, then by resorting to subsequences of $< k_i >$ if necessary, we assume that, for some points $u$ and $v$ in $[z, c] \cup [a, b]$, $\lim_{i \to \iy} g^{k_i-1}(x) = u$ and $\lim_{i \to \iy} g^{k_i-1}(y) = v$.  Since $z$ and $b$ are the only two points in $[z, c] \cup [a, b]$ which are mapped to $z$ under $g$, we see that $u = z$ or $u = b$.  Consequently, $f^j(v) \ne u$ for all $1 \le j \le m-1$.  Let $\dl = \min \{ \, |b-z|, |f^j(v)-u| \, \big| \, 1 \le j \le m-1 \, \}$.  Then, $\lim_{n \to \iy} |f^n(f^s(x))-f^n(f^t(y))| \ge$ $\lim_{i \to \iy} |f^{m(k_i-1)}(x) - f^{m(k_i-1)+t-s}(y)| = \lim_{i \to \iy} |g^{k_i-1}(x) - f^{t-s}(g^{k_i-1}(y))| = |u - f^{t-s}(v)| \ge \dl$.  The rest is easy and omitted.
\hfill\sq

\noindent
{\bf{Remark 1.}}
Let $f(x)$ be a continuous map from $[-1, 1]$ onto itself defined by $f(x) = 2x+2$ for $-1 \le x \le -1/2$; $f(x) = -2x$ for $-1/2 \le x \le 0$; and $f(x) = -x$ for $0 \le x \le 1$.  Then, by Theorem 3, $f$ has uncountable invariant scrambled sets although $f$ is not turbulent.

\section{Maps with periodic points of periods not a power of 2.}
If there are closed subintervals $J_0$, $J_1$, $\cdots$, $J_{n-1}, J_n$ of $I$ with $J_n = J_0$ such that $f(J_i) \supp J_{i+1}$ for $i = 0, 1, \cdots, n-1$, then we say that $J_0J_1 \cdots J_{n-1}J_0$ is a {\it{cycle of length}} $n$.  It is well known that if $J_0J_1\cdots J_{n-1}J_0$ is a cycle of length $n$ then there exists a periodic point $y$ of $f$ such that $f^i(y)$ belongs to $J_i$ for $i = 0, 1, \cdots, n-1$ and $f^n(y) = y$.  The following result is well known [{\bf 1}, p. 26].  Here we give a different proof which is interesting in its own right. 

\noindent
{\bf Theorem 4.}
{ \it If $f$ has a periodic point of odd period $m \ge 3$ and no periodic points of odd period $k$ with $1 < k < m$, then any periodic orbit of $f$ with least period $m$ must be a \v Stefan orbit (see {\rm [{\bf 1}, p. 11]} for the definition).  Consequently, $f^n$ is strictly turbulent for $n = 2$ and every $n \ge m$.}

\noindent
{\it Proof.}
Let $P = \{ x_i \big| 1 \le i \le m \}$, with $x_1 < x_2 < \cdots < x_m$, be a periodic orbit of $f$ with least period $m$.  If $m = 3$, the proof is easy.  So, suppose $m > 3$.  Let $x_s = \max \{ x \in P \big| x < f(x) \}$.  Then $f(x_{s+1}) \le x_s$ and $x_{s+1} \le f(x_s)$ and so $f$ has a fixed point $z$ in $[x_s, x_{s+1}]$.  Since $m$ is odd, for some integer $1 \le t \le m-1$ and $t \ne s$, $f(x_t)$ and $f(x_{t+1})$ lie on opposite sides of $z$.  For simplicity, we assume that $x_t < x_s$.  If $x_{s+1} \le x_t$, the proof is similar.  

If, for some $2 \le q \le m-3$, $f^q(x_s) \le x_t$, then let $J_i = [f^i(x_s):z]$ for all $0 \le i \le q-1$ and consider the cycle $J_0J_1 \cdots J_{q-1}[x_t, x_{t+1}]J_0 \cdots J_0$ of odd length $m-2$.  If, for some $1 \le r \le m-3$, $f^r(x_s)$ and $f^{r+1}(x_s)$ lie on same side of $z$, then since $f$ maps each of $x_s$ and $x_{s+1}$ to the other side of $z$, $f^r(x_s) \notin \{ \, x_s, x_{s+1}\, \}$.  Let $J_i = [f^i(x_s):z]$ for all $0 \le i \le r-1$ and $J_r = [f^r(x_s), x_s]$ if $f^r(x_s)< x_s$ and $J_r = [x_{s+1}, f^r(x_s)]$ if $x_{s+1} < f^r(x_s)$ and consider the cycle $J_0J_1 \cdots J_rJ_0 \cdots J_0$ of odd length $m-2$.  In either case, we obtain a periodic point of $f$ of odd period strictly between 1 and $m$.  This is a contradiction.  So, $f^{m-1}(x_s) = x_t < f^i(x_s) < z < f^j(x_s)$ for all even $i$ in $[0, m-3]$ and all odd $j$ in $[1, m-2]$.    

For all $0 \le i \le m-1$, let $J_i = [z:f^i(x_s)]$.  If, for some $1 \le j < k \le m-1$ with $k-j$ even, $[z : f^k(x_s)] \sub [z : f^j(x_s)]$, then by considering the cycle $J_0J_1 \cdots J_{j-1}J_kJ_{k+1} \cdots$ $J_{m-2}$ $[x_t, x_{t+1}]J_0$ of odd length $m - k + j < m$, we obtain a periodic point of $f$ of odd period strictly between 1 and $m$.  This contradicts the assumption.  So, we must have     
$$
f^{m-1}(x_s) < \cdots < f^4(x_s) < f^2(x_s) < x_s < f(x_s) < f^3(x_s) < \cdots < f^{m-2}(x_s).
$$
\noindent
That is, the orbit of $x_s$ is a \v stefan orbit and $x_s < z < x_{s+1} = f(x_s)$.  Let $w$ be a point in $[f^{m-4}(x_s)$, $f^{m-2}(x_s)]$ such that $f^2(w) = z$.  Let $z < u < f^{m-4}(x_s) < v < w$ be points such that $f^2(u) = f^2(v) = w$.  Let $J_0 = [z, u]$ and $J_1 = [v, w]$.  Then $J_0 \cap J_1 = \emp$ and $f^2(J_0) \cap f^2(J_1) \supp J_0 \cup J_1$.  So, $f^2$ is strictly turbulent.  On the other hand, by direct computations, $f^{m-2}([x_s, f(x_s)])$ $\cap f^{m-1}([f^{m-1}(x_s), f^{m-3}(x_s)])$ $ \supp [\min P, \max P]$.  Thus, $f^n$ is strictly turbulent for any $n \ge m-1$.
\hfill\sq

By Sharkovsky's theorem {\bf{\cite{du2}}} and Theorems 1 $\&$ 4, we obtain the following result.

\noindent
{\bf Theorem 5.}
{\it If $f$ has a periodic point of least period $2^k \cdot m$, where $m \ge 3$ is odd and $k \ge 0$,  then for $n = 2$ and every $n \ge m$ there exist a number $\dl_n > 0$ and an uncountable $\dl_n$-scrambled set $S_n$ of $f$ in the recurrent set of $f$ which is invariant under $f^{2^k \cdot n}$ such that $\bigcap_{x \in S_n} \omega(x, f)$ contains $S_n$ and contains periodic points of $f$ of arbitrarily large periods.}

\noindent
{\bf{Remark 2.}}
For every positive integer $i$, let $f_i$ be the continuous map from $[1, 2i+1]$ onto itself with the following eight properties: (1) $f_i(1) = i+1$; (2) $f_i(2) = 2i+1$; (3) $f_i(i+1) = i+2$; (4) $f_i(i+2) = i$; (5) $f_i(i+1+1/3) = i+1+2/3$; (6) $f_i(i+1+2/3) = i+1+1/3$; (7) $f_i(2i+1) = 1$; and (8) $f_i$ is linear on each of the intervals: $[1, 2]$, $[2, i+1]$, $[i+1, i+1+1/3]$, $[i+1+1/3, i+1+2/3]$, $[i+1+2/3, i+2]$ and $[i+2, 2i+1]$.  Then $|f_i(x)-x| \ge 1/3$ for any point $x$ in $[1, i+1+1/3]$ and in $[i+1+2/3, 2i+1]$ and $f_i^2(x) = x$ for all $x$ in $[i+1+1/3, i+1+2/3]$.  So, $f_i$ cannot have any non-empty invariant scrambled set while $\{1, 2, \cdots, 2i+1\}$ is a \v Stefan orbit of $f_i$ with least period $2i+1$.  By successive applications to $f_i$ of the square root construction as described in section 1, we see that the result in Theorem 5 is best possible in the sense that there exist continuous interval maps which have periodic points of least period $2^k \cdot m$ with odd $m \ge 3$, but no uncountable scrambled sets which are invariant under $f^{2^k}$.  

\section{The\ shift\ map\ on\ $\Si_2$\ has\ dense\ uncountable\ invariant\ 1-scrambled\ sets.} 
Let $\al = \al_0\al_1\al_2 \cdots$ be a transitive point of $\Si_2$.  For any $\be = \be_0\be_1\be_2 \cdots$ in $\Si_2$, let 

\noindent
$\tau_\be = \, \al_0\be_0 \, (01) \, \al_0\al_1(\be_0)^2(\be_1)^2 \, (01)^2(0011)^2\al_0\al_1\al_2(\be_0)^3(\be_1)^3(\be_2)^3 \, (01)^3(0011)^3(000111)^3$

\quad $\al_0\al_1\al_2\al_3(\be_0)^4(\be_1)^4(\be_2)^4(\be_3)^4 \, (01)^4(0011)^4(000111)^4(00001111)^4 \cdots,$

\noindent
where $(\be_0)^2 = \be_0\be_0$, $(01)^3 = 01 \, 01 \, 01$, $(0011)^4 = 0011 \, 0011 \, 0011 \, 0011$, etc.  Let $S = \{ \, \si^k(\tau_\be) \, \mid \, \be \in \Si_2, k \ge 0 \, \}$.   Then we easily obtain the following result (note that, in $\Si_2$, the longest distance between any two points is 1).

\noindent
{\bf Theorem 6.}
{\it The set $S$ is a dense uncountable invariant $1$-scrambled set for $\si$ which consists of transitive points.}


\noindent
\begin{thebibliography}{99}
\bibitem{bc} L. S. Block and W. A. Coppel, \it Dynamics in One Dimension, \rm Lecture Notes in Math., no. 1513, Springer-Verlag, Berlin, 1992.

\bibitem{de} R. L. Devaney, \it An introduction to chaotic dynamical systems, \rm 2nd edition, Addison-Wesley, Redwood City, CA., 1989. 

\bibitem{du1} B.-S. Du, Every chaotic interval map has a scrambled set in the recurrent set, \it Bull. Austral. Math. Soc. \rm {\bf 39} (1989), 259-264.

\bibitem{du2} B.-S. Du, A simple proof of Sharkovsky's theorem, \it Amer. Math. Monthly \rm {\bf 111} (2004), 595-599.

\bibitem{ev} M. J. Evens, P. D. Humke, C.-M. Lee and R. J. O'malley, Characterizations of turbulent one-dimensional mappings via $\omega$-limit sets, \it Trans. Amer. Math. Soc. \rm {\bf 326} (1991), 261-280; Corrigendum: \it Trans. Amer. Math. Soc. \rm {\bf 333} (1992), 939-940.

\bibitem{ly} T.-Y. Li and J. A. Yorke, Period three implies chaos, \it Amer. Math. Monthly \rm {\bf 82} (1975), 985-992.
\end{thebibliography}
\end{document}